\documentclass{article}
\usepackage{amsmath,amsthm,amsfonts,amssymb}
\usepackage[all]{xypic}

\theoremstyle{plain}
\newtheorem{prop}{Proposition}
\newtheorem{thm}{Theorem}
\newtheorem{cor}{Corollary}
\newtheorem{lem}{Lemma}

\theoremstyle{definition}
\newtheorem{remark}{Remark}
\newtheorem{example}{Example}
\newtheorem{defi}{Definition}

\newcommand{\C}[1]{{\cal#1}} % Calligraphic

\newcommand{\vsp}{\vspace{0.05in}}
\newcommand{\bea}{\begin{eqnarray}}
\newcommand{\eea}{\end{eqnarray}}
\newcommand{\noi}{\noindent}
\newcommand{\nn}{\nonumber}

\begin{document}

\title{A Survey of Huebschmann and Stasheff's Paper: 
Formal solution of the master equation via HPT and deformation theory} 
\author{F\"{u}sun Akman and Lucian M.~Ionescu}

\maketitle 

%\tableofcontents

\section{Introduction}

These notes, based on the paper~\cite{HS} by Huebschmann and Stasheff, were prepared for a series of talks at Illinois State University with the intention of applying Homological Perturbation Theory (HPT) to the construction of derived brackets~\cite{Kosmann,Vor}, and eventually writing Part~II of the paper~\cite{AI}.

Derived brackets are obtained by deforming the initial bracket via a derivation of the bracket. In~\cite{AIPapa} it was demonstrated that such deformations correspond to solutions of the Maurer-Cartan equation, and the role of an ``almost contraction'' was noted. This technique (see also~\cite{I-comb}) is very similar to the iterative procedure of~\cite{HS} for finding the most general solution of the Maurer-Cartan equation, i.e.~the deformation of a given structure in a prescribed direction.

The present article, besides providing additional details of the condensed article \cite{HS}, forms a theoretical background for understanding and generalizing the current techniques that give rise to derived brackets. The generalization, which will be the subject matter of \cite{AI-2}, will be achieved by using Stasheff and Huebschmann's universal solution. A second application of the universal solution will be in deformation quantization and will help us find the coefficients of star products in a combinatorial manner,
rather than as a byproduct of string theory which underlies the original solution given by Kontsevich \cite{Kon}.

HPT is often used to replace given chain complexes by homotopic, smaller, and more readily computable chain complexes (to explore ``small'' or ``minimal'' models). This method may prove to be more efficient than ``spectral sequences'' in computing (co)homology. One useful tool in HPT is
\begin{lem}[Basic Perturbation Lemma (BPL)] Given a contraction of $N$ onto $M$ and a perturbation $\partial$ of $d_N$, under suitable conditions there exists a perturbation $d_{\partial}$ of $d_M$ such that $H(M,d_M+d_{\partial})=H(N,d_N+\partial)$.
\end{lem}
The main question is: under what conditions does the BPL allow the preservation of the data structures (DGA's, DG coalgebras, DGLA's etc.)? (We will use the self-explanatory abbreviations such as DG for ``differential graded'', DGA for ``differential graded (not necessarily associative) algebra'', and DGLA for ``differential graded Lie algebra''.) 

Another prominent idea is that of a ``(universal) twisting cochain'' 
as a solution of the ``master equation'':
\begin{prop} Given a contraction of $N$ onto $M$ and a twisting cochain $N\rightarrow A$ ($A$ some DGA), there exists a unique twisting cochain $M\rightarrow A$ that factors through the given one and which can be constructed inductively.
\end{prop}
%Such universal solutions of the master equation on the contracted space 
%can be obtained inductively from the solutions on the big space. (The Proposition says that!)
The explicit formulas are reminiscent of the Kuranishi map \cite{Manetti} (p.17), and the
relationship will be investigated elsewhere.

Note: we will assume that the ground ring is a field $F$ of characteristic zero. 
We will denote the end of an example with the symbol $\Diamond$ and the end of a proof by $\Box$.%Results can be mostly extended to the case of a ground ring containing the rational numbers as a subring. 

% **********************************************************************
\section{Perturbations of (co)differentials}

\subsection{Derivations of the tensor algebra}

For any vector space $V$ over $F$ we have the isomorphism $\mbox{Der}(TV)\cong \mbox{Hom}(V,TV)$ where $TV$ denotes the (augmented) tensor algebra on $V$. Namely, every linear map $f$ from $V$ into $TV$ extends uniquely into a derivation of the algebra $TV$ via the formula
\[ \hat{f}(v_1\otimes\cdots v_n)=\sum_{i=1}^nv_1\otimes\cdots\otimes f(v_i)\otimes\cdots v_n. \]
Equivalently, every derivation of $TV$ is determined by its restriction to $V$.

\subsection{Coderivations of the tensor coalgebra}

Similarly, we have the isomorphism $\mbox{Coder}(T^cV)\cong \mbox{Hom}(T^cV,V)$ where $T^cV$ is the (coaugmented) coassociative tensor coalgebra of $V$, with counit
$\eta:T^cV\rightarrow F$ (projection onto $F$), and comultiplication
\[ \Delta (v_1\otimes\cdots\otimes v_n)=\sum_{i=0}^n(v_1\otimes\cdots\otimes v_i)\otimes(v_{i+1}\otimes\cdots\otimes v_n).\]
Every linear map $f=f_1+f_2+\cdots +f_n+\cdots :T^cV\rightarrow V$ (where $f_i:V^{\otimes i}\rightarrow V$) factors uniquely through a coderivation $\hat{f}$ of $T^cV$ defined via the formula
\begin{eqnarray}
\hat{f}(v_1\otimes\cdots\otimes v_n)&=&\sum_{i=1}^nv_1\otimes\cdots\otimes f_1(v_i)\otimes\cdots v_n \nonumber\\
&&+\sum_{i=1}^{n-1}v_1\otimes\cdots\otimes f_2(v_i\otimes v_{i+1})\otimes\cdots v_n \nonumber\\
&&\vdots\nonumber\\
&&+f_n(v_1\otimes\cdots\otimes v_n).\nonumber\end{eqnarray}
That is, each coderivation on $T^cV$ is determined by itself followed by the projection onto $V$. Recall that the condition for $\hat{f}$ to be a coderivation can be written as $ \Delta \hat{f}=(1\otimes \hat{f}+\hat{f}\otimes 1)\Delta .$

\subsection{Coderivations of the symmetric coalgebra}\label{twothree}

Let us consider the cofree cocommutative counital coassociative algebra $ST^cV$ on the vector space $V$ as a subspace of $T^cV$. The symmetric group $\Sigma_n$ acts on the left on $V^{\otimes n}$ via $\sigma(v_1\otimes\cdots\otimes v_n)=v_{\sigma^{-1}(1)}\otimes\cdots\otimes v_{\sigma^{-1}(n)}$. Then 
\[ ST^cV=\bigoplus_{n\geq 0}(V^{\otimes n})^{\Sigma_n}\]
is the space of invariants of this action. The action is compatible with the coproduct
on $ST^cV$, so $ST^cV$ is a subcoalgebra of $T^cV$ which is cocommutative. Note that $ST(V)$ is not a subalgebra with respect to the tensor multiplication in $T(V)$; the product has to be symmetrized so that it projects back onto this subspace (reminiscent of what T.~Voronov does with derived brackets). The projection (symmetrization) map $P:T^cV\rightarrow ST^cV$ is given by
\[ P(v_1\otimes\cdots\otimes v_n)=\frac{1}{n!}\sum_{\sigma\in\Sigma_n}\sigma(v_1\otimes\cdots\otimes v_n).\]
This is not a coalgebra map, but is a retraction of the canonical inclusion $ST^cV\hookrightarrow T^cV$. Now a coderivation $D:T^cV\rightarrow T^cV$ induces a coderivation $D_S:ST^cV\rightarrow ST^cV$ by the composition $ST^cV\hookrightarrow T^cV\stackrel{f}{\rightarrow} T^cV\stackrel{P}{\rightarrow} ST^cV$. In particular, a coderivation of $T^cV$ induces one of $\Lambda^cV=ST^c[sV]$ where $V$ is thought of as living in degree zero: we introduce the {\it graded} symmetric coalgebra below. Once again, coderivations of $ST^cV$ are determined by their projections onto $ST^c_1V$; a map $f=f_1+f_2+\cdots:ST^c(V)\rightarrow V$ determines a coderivation $\hat{f}$ as in the tensor coalgebra case.
\vsp

In the remaining part of this survey, we choose to identify $ST^cV$ with the abstract symmetric coalgebra $S^cV$ under the isomorphism
\[ v_1\cdots v_n\mapsto P(v_{1}\otimes\cdots\otimes v_{n}).\]
The coproduct in $S^cV$ is given by
\[ \Delta(v_1\cdots v_n)=\sum_{i=0}^n\sum_{\sigma\in\Sigma_{i,n-i}}v_{\sigma(1)}\cdots v_{\sigma(i)}\otimes v_{\sigma(i+1)}\cdots v_{\sigma(n)}.\]

\subsection{DGLA's and perturbations of the codifferential}

\begin{defi} For any chain complex $(X,d)$, and odd $\partial$, with $(d+\partial)^2=0$, we say that $\partial$ is a {\it perturbation} of the differential $d$. We call $d+\partial$ the {\it perturbed differential}. Equivalently, we have $[d,\partial]+\partial\partial=0$ in $\mbox{End}(X)$. If $\partial$ is also compatible with an existing coalgebra structure on $X$, we say that it is a {\it coalgebra perturbation}. 
\end{defi}
Let $(g,d)$ be a graded chain complex ($d$ lowers degrees) with a bracket $[\,,\,]$ that is skew-symmetric (not necessarily Leibniz or a chain map). Consider the differential graded symmetric coalgebra $S^c[sg]$, the differential $d$ being induced by that on $g$. Also let $\partial$ be the coderivation on $S^c[sg]$ of degree $-1$ induced by the bracket.
\begin{prop} The bracket $[\,,\,]$ turns $(g,d)$ into a DGLA if and only if $\partial$ is a coalgebra perturbation of $d$. Also, any DGLA structure on $g$ is determined by the coalgebra perturbation induced from the bracket.
\end{prop}
When $g$ is an ordinary (degree-zero) Lie algebra over a field, $S^c[sg]=\Lambda^cg$ with differential $\partial$ corresponding to the bracket is the ordinary Koszul or Chevalley-Eilenberg complex computing the homology of $g$ with coefficients in the field.

\subsection{Strongly homotopy Lie algebras}

\begin{defi} Let $(g,d)$ be a chain complex and let $d$ also denote the codifferential in $S^c[sg]$ induced by $d$. A {\it strongly homotopy Lie} (sh-Lie, or $L_{\infty}$) structure on $g$ is a perturbation $\partial=\partial_2+\cdots +\partial_n+\cdots$ of $d$, i.e.~an odd coderivation satisfying $[d,\partial]+ \partial\partial =0$ and $\partial\eta =0$ (recall that $\eta$ is the counit) so that the sum $d+\partial$ endows $S^c[sg]$ with a new coaugmented DG coalgebra structure.\end{defi}
The corresponding mega-map $\ell_2+\cdots +\ell_n+\cdots$ from $S^c(sg)$ to $g$ extends the differential $\ell_1=d:sg\rightarrow g$, and the lower identities satisfied by 
\[ \ell=\ell_1+\ell_2+\cdots +\ell_n+\cdots\]
read as follows:
\bea &&\ell_1^2=0\nn\\
&&\ell_1(\ell_2(a,b))\pm\ell_2(\ell_1(a),b)\pm\ell_2(\ell_1(b),a)=0\nn\\
&&\ell_1(\ell_3(a,b,c))\pm \ell_3(\ell_1(a),b,c)\pm\ell_3(\ell_1(b),a,c)\pm\ell_3(\ell_1(c),a,b)\nn\\
&&\pm\ell_2(\ell_2(a,b),c)\pm\ell_2(\ell_2(a,c),b)\pm\ell_2(\ell_2(b,c),a)=0.\nn\eea
An sh-Lie morphism between two sh-Lie (or DGL) algebras $(g,d+\cdots)$ and $(g',d'+\cdots)$ is a collection of chain maps $F_n:S^c_n[sg]\rightarrow S^c_n[sg']$, satisfying $\Delta' F(u)=(F\otimes F)(\Delta u)$. Then $F$ is uniquely determined by its projection onto $sg'$, that is, we may assume $F_n:S^c_n[sg]\rightarrow sg'$. 
\begin{defi} A {\it quasi-isomorphism} $F$ between sh-Lie algebras $g,g'$ is an sh-Lie morphism such that $F_1:sg\rightarrow sg'$ induces an isomorphism between $H(g,d)$ and $H(g',d')$. 
\end{defi}
\begin{remark} Quasi-isomorphisms between DGLA's are especially important in deformation theory. Such a map gives a one-to-one correspondence between moduli spaces of solutions to MC equations in $\hbar g[[\hbar]]$ and $\hbar g'[[\hbar]]$ (see~\cite{Chen}): given a quasi-isomorphism $F:S^c[sg]\rightarrow sg'$, we define $\tilde{F}:\hbar g[[\hbar]]\rightarrow \hbar g'[[\hbar]]$ by
\[ \tilde{F}(r)=\sum_{n=1}^{\infty}\frac{1}{n!}F_n(r,\dots,r)\]
(also see~ \cite{Dol}).
\end{remark}

\subsection{The Hochschild chain complex and DGA's}

Let $(A,\mu)$ be a unital associative algebra (possibly graded), and $T^c[sA]$ denote the tensor coalgebra on the suspension of $A$. We recall that 
\[ \mbox{Coder}(T^c[sA])\cong \mbox{Hom}(T^c[sA],A). \]
In particular, the associative bilinear multiplication $\mu\in\mbox{Hom}(T^c[sA],A)$ corresponds to a square-zero coderivation $\partial:T^c[sA]\rightarrow T^c[sA]$ defined by
\bea && \partial(a_1\otimes \cdots\otimes a_n)\nn\\ &=&\sum_{i=1}^{n-1}(-1)^{i+1}(a_1\otimes \cdots\otimes\mu(a_i\otimes a_{i+1})\otimes\cdots\otimes a_n)\nn\\ &&+(-1)^{n+1}(\mu(a_n\otimes a_1)\otimes a_2\otimes \cdots\otimes a_{n-1}) .   \nn\eea
The condition that $\partial$ is a codifferential is equivalent to the associativity condition $m\circ m=0$ where $\circ$ is the Gerstenhaber composition on multilinear maps (a right pre-Lie map). The complex $(\mbox{Hom}(T^c[sA],A),\partial)$ is known as the {\it Hochschild chain complex}. 

Now let $(A,\mu,d)$ be a DGA. Then $d+\mu\in \mbox{Hom}(T^c[sA],A)$ corresponds to a perturbed codifferential $d+\partial$ satisfying $(d+\partial)^2=0$, which is equivalent to the identities $d^2=0$ and $[d,\partial]+\partial\partial =0$. The latter can also be split into $[d,\partial]=0$ and $\partial\partial =0$. 
\begin{prop} The multiplication $\mu$ turns $(A,d)$ into a DGA if and only if $\partial$ is a coalgebra perturbation of $d$. Also, any DGA structure on $A$ is determined by the coalgebra perturbation induced from $\mu$.
\end{prop}

%This section needs to be better understood. What is the relationship between the homology and cohomology differentials?\vsp
%
%[According to Loday~\cite{Lo}, Connes discovered in 1981 the following phenomenon: the Hochschild boundary map is still well-defined when we factor out $A^{\otimes n}$ by the action of $C_n$. That is, it is possible to define {\it cyclic homology} where the chain complex is $C^c[sA]$ instead of $T^c[sA]$.]

\subsection {Strongly homotopy associative algebras}
 
\begin{defi} Let $(A,d)$ be a chain complex and let $d$ also denote the codifferential in $T^c[sA]$ induced by $d$. A {\it strongly homotopy associative} (or $A_{\infty}$) structure on $A$ is a perturbation $\partial=\partial_2+\cdots +\partial_n+\cdots$ of $d$, i.e.~an odd coderivation satisfying $[d,\partial]+ \partial\partial =0$ and $\partial\eta =0$ so that the sum $d+\partial$ endows $T^c[sA]$ with a new coaugmented DG coalgebra structure.\end{defi}
The corresponding mega-map $m_2+\cdots +m_n+\cdots$ from $T^c[sA]$ to $A$ extends the differential $m_1=d:sA\rightarrow A$, and the lower identities satisfied by 
\[ m=m_1+m_2+\cdots +m_n+\cdots\]
read as follows:
\bea &&m_1^2=0\nn\\
&&m_1(m_2(a,b))\pm m_2(m_1(a),b)\pm m_2(m_1(b),a)=0\nn\\
&&m_1(m_3(a,b,c))\pm m_3(m_1(a),b,c)\pm m_3(a,m_1(b),c)\pm m_3(a,b,m_1(c))=0.\nn\eea
The mega-identity is $m\circ m=0$, sometimes written in the braces notation $\{ m\}\{ m\} =0$.

\section{Master equation}

If $(A,d)$ is a differential graded associative algebra (DGA), then the equation
\begin{equation} d\tau =\tau\tau \label{one}\end{equation}
is called the {\it Master Equation (ME)} (or {\it Maurer-Cartan equation (MCE)}, etc.). Similarly if $(g,d)$ is a DGLA, then the equation
\begin{equation} d\tau =\frac{1}{2}[\tau,\tau]\label{two}\end{equation}
is also called the Master Equation. Sometimes the sign convention is
\begin{equation}  d\tau +\tau\tau=0\label{three}\end{equation}
or
\begin{equation}  d\tau +\frac{1}{2}[\tau,\tau]=0.\label{four}\end{equation}
Clearly any solution of such an equation must be an odd element of the algebra. Moreover, in case $A$ is the graded universal enveloping algebra of $g$, or $g$ is the Lie algebra obtained from $A$ by the usual bracket, then solutions of the DGLA master equation are also solutions of the DGA master equation.
\begin{remark} If $\tau$ is a solution of Eq.~(\ref{two}) or Eq.~(\ref{four}) in a DGLA $g$, then the odd derivation $d_{\tau}=d-\mbox{ad}_{\tau}$ or $d_{\tau}=d+\mbox{ad}_{\tau}$ respectively defines a new DGLA structure on $g$ with respect to the old bracket. \end{remark}
\begin{example} If $g$ is a DGLA or $L_{\infty}$ algebra, then the corresponding coalgebra perturbation $\partial$ in $\mbox{Coder}(S^c(sg))$ is a solution of the ME in the DGA $\mbox{End}(S^c(sg))$, where the differential is 
$D=\mbox{ad}_d$. $\diamond$\end{example}
\begin{example} Gauge Theory: Let $\xi$ be a principal bundle with structure group $G$ and Lie algebra $g$. There is a graded Lie algebra structure on the ad$(\xi)$-valued de Rham forms induced by $g$. Given a connection $A$ and an ad$(\xi)$-valued 1-form $\eta$, the sum $A+\eta$ is again a connection, and its curvature is
\[ F_{A+\eta}=F_A+d_A\eta +\frac{1}{2}[\eta,\eta].\]
In particular, $F_A=F_{A+\eta}$ if and only if
\[ d_A\eta +\frac{1}{2}[\eta,\eta]=0\]
(the Maurer-Cartan equation). Here $d_A$ is the covariant derivative of the connection $A$. When $A$ is a flat connection (zero curvature) then there exists a DGLA structure on the ad($\xi$)-valued differential forms ($d_A^2=0$) and $F_{A+\eta}$ is also flat iff the MCE is satisfied (then the covariant derivative for $A+\eta$ is $d_{\tau}$). $\diamond$
\end{example}

\section{Twisting cochain}

The notion of a twisting cochain generalizes that of a connection in differential geometry. 

\subsection{Differential on Hom}

If $(C,d_C)$ and $(A,d_A)$ are chain complexes, the following differential $D$ makes $\mbox{Hom}(C,A)$ into a chain complex: $D\phi =d_A\phi\pm \phi d_C$.

\subsection{Cup product and cup bracket}

\begin{prop} For any differential graded coalgebra $C$ and a differential graded associative algebra (DGA) $A$, the chain complex $(\mbox{Hom}(C,A),D)$ becomes a DGA via the {\it cup (convolution) product} $a\smile b$ defined by the composition
\[ C\stackrel{\Delta}{\longrightarrow} C\otimes C\stackrel{a\otimes b}{\longrightarrow} A\otimes A\stackrel{\mu}{\longrightarrow} A.\]
The coaugmentation and augmentation maps $\eta$ and $\epsilon$ on $C$ and $A$ respectively define an augmentation map on $(\mbox{Hom}(C,A),D)$.
\end{prop}
\begin{prop} For any differential graded coalgebra $C$ and a DGLA $g$, the chain complex $(\mbox{Hom}(C,g),D)$ becomes a DGLA via the {\it cup bracket} $[a,b]$ defined by the composition
\[ C\stackrel{\Delta}{\longrightarrow}C\otimes C\stackrel{a\otimes b}{\longrightarrow} g\otimes g\stackrel{[\,,\,]}{\longrightarrow} g.\]
\end{prop}
\begin{example} If $g$ is a Lie algebra, then the cup bracket on $\mbox{Hom}(S^c[sg],g)$ is defined as above. For example, if $\tau$ and $\kappa$ are maps $S^c_1[g]\rightarrow g$, then $[\tau,\kappa]$ may be nonzero only on $S^c_2[sg]$. In this case, we compute
\[ \Delta(xy)=1\otimes xy+x\otimes y+y\otimes x+xy\otimes 1\;\;(x,y\in sg) ,\]
and 
\begin{equation} [\tau,\kappa](xy)=[\tau(x),\kappa(y)]+[\tau(y),\kappa(x)].\label{cupbra}
\end{equation}
$\diamond$
\end{example}
\begin{example} %{\bf (May have to add, or switch to, the Hochschild homology complex. See pp. 37,51 of Loday.)} 
The Hochschild complex of an associative algebra $(A,\mu)$ (where $\mu^2=0$): Let
\[ C^{\bullet}(A)=\mbox{Hom}(T^cA,A)=  \mbox{Hom}(\bigoplus_{n\geq 0}A^{\otimes n},A),\]
with differential $D=\mbox{ad}_{\mu} \in\mbox{Der}(C^{\bullet}(A))$. The cup product $x\smile y=\{ \mu\}\{ x, y\}$ is the composition
\[ T^cA\stackrel{\Delta}{\longrightarrow}T^cA\otimes T^cA\stackrel{x\otimes y}{\longrightarrow}A\otimes A \stackrel{\mu}{\longrightarrow} A;\]
if $x$ is an $n$-linear map and $y$ is an $m$-linear map, then 
\[ (x\smile y)(a_1\otimes\cdots\otimes a_{n+m})=x(a_1\otimes\cdots\otimes a_{n})\cdot y(a_{n+1}\otimes\cdots\otimes a_{n+m}).\]
$\diamond$\end{example}
\begin{remark} The differential $D$ above is an inner derivation and not derived from differentials on $A$ and $T^cA$. Still, it is a derivation of the cup product. %These maps and more have been studied by Gerstenhaber and a few thousand others. The structure of a $G$-algebra I mentioned earlier appears on $\mbox{End}(C^{\bullet}(A))$ with $M(1)$ as the differential, $M(2)$ as the cup product, $M(1|1)$ as the substitution operator (right pre-Lie), and higher substitution operators $M(1|n)$. The anti-symmetrization of $M(1|1)$ becomes a Gerstenhaber bracket in the $M(1)$-cohomology.
\end{remark}

\subsection{Twisting cochain}

\begin{defi} Given a coaugmented DG coalgebra $C$ and an augmented DGA $A$, a {\it twisting cochain} is a homogeneous morphism $t:C\rightarrow A$ of degree $-1$ such that $\epsilon\tau=0$ and $\tau\eta=0$ , and which satisfies $Dt=t\smile t$.
\end{defi}
In other words, a twisting cochain is a solution of the master equation on $\mbox{Hom}(C,A)$ with the usual differential $D$ induced from those of $C$, $A$ and the product is the cup product. 
%\noi {\bf Example.} If the funny derivation is allowed, then solutions of the master equation on the Hochschild complex are odd, linear maps $T$ on $A$ with
%\[ T(ab)\pm T(a)b\pm aT(b)\pm T(a)T(b)=0.\]
%Any significance? {\bf Notes:} (Lucian) 1.~The differential $D$ is a perturbation of the zero differential. 2.~This looks like the Rota-Baxter Equation. $\diamond$\vsp
\begin{defi} Given a DG cocommutative coalgebra $C$ and a DGLA $g$, a {\it Lie algebra twisting cochain} $t:C\rightarrow g$ is a homogeneous map of degree $-1$ whose composition with the coaugmentation map is zero, and which satisfies
\begin{equation} Dt=\frac{1}{2}[t,t]\label{twico}\end{equation}
($[\,,\,]$ being the cup bracket).
\end{defi}
Recall that a DGLA structure on a graded chain complex $(g,d)$ is given by a perturbation $\partial$ of the corresponding codifferential on $S^c[sg]$. Moreover, the piece $d\partial+\partial d=0$ of $(d+\partial)^2=0$ says that the bracket is a chain map and the piece $\partial^2=0$ says that the bracket satisfies the Jacobi identity. Let us denote the symmetric coalgebra with the codifferential $\partial$ by $S_{[\, ,\,]}^c[sg]$. Quillen's notation $C[g]$ for the same DG coalgebra reminds us that this is the Koszul or Chevalley-Eilenberg complex that computes the homology of $g$ without any regard for the additional differential $d$ on $g$.  
\begin{example} For any DGLA $g$, its {\it universal Lie algebra twisting cochain}
\[ \tau_g:S^c_{[\, ,\,]}[sg]\rightarrow g\]
is given by 
\[ \begin{array}{ll} \tau_g(sx)=x & \mbox{for $x\in g$}\\ \tau_g(y)=0 & \mbox{for $y\in S^c_k[sg]$, $k\neq 1$.}\end{array}\]
That is, an element with tensor degree one goes to its desuspension and everything else goes to zero. Clearly, the composition $\tau_g\eta$ is zero, as $\tau_g=0$ on constants. Next, we show that $\tau_g$ satisfies the equation~(\ref{twico}), but we note that in this construction the differential on $g$ itself is taken to be zero. On the left-hand side, we have
\[ D\tau_g(x_1\wedge \cdots \wedge x_n)=\tau_g\partial (x_1\wedge \cdots \wedge x_n)\]
which is zero if $n\neq 2$ and is equal to $[x_1,x_2]$ if $n=2$.
Meanwhile, on the right-hand side, we have
\[ \frac{1}{2}[\tau_g,\tau_g](x_1\wedge \cdots \wedge x_n)\]
which is zero if $n\neq 2$ and is equal to 
\[ \frac{1}{2}\{ [\tau_g(x_1),\tau_g(x_2)]\pm [\tau_g(x_2),\tau_g(x_1)]\} =[x_1,x_2]\]
if $n=2$. $\diamond$
\end{example}
\begin{remark} The universal property of the universal LA twisting cochain is that every Lie algebra twisting cochain factors through this one: that is, whenever $C$ is a coalgebra and $\tau:C\rightarrow g$ is a twisting cochain, then $\tau_g\circ c(\tau)=\tau$ where $c(\tau):C\rightarrow S^c_{[\, ,\,]}[sg]$ is the unique coalgebra map induced by $\tau$.
\end{remark}
Using HPT, we will construct formal solutions 
\[ \tau\in \mbox{Hom}(S^c_D[sH(g)],g)\]
of the master equation. Once we make the choice of a contraction, we will obtain explicit inductive formulas for $D$ and $\tau$.

\section{Homological perturbation theory (HPT)}

``HPT is concerned with transferring various kinds of algebraic structure through a homotopy equivalence''. Also: ``HPT is a set of techniques for the transference of structures from one object to another up to homotopy'' (Real~\cite{Re}).

\subsection{Contraction}

\begin{defi} Let $(M,d_M)$ and $(N,d_N)$ be chain complexes, $\pi:N\rightarrow M$ and $\nabla:M\rightarrow N$ be chain maps, and $h\in\mbox{End}(N)$ be a morphism (possibly preserving some extra structure) of degree~1. Then a {\it contraction}
\begin{equation} (   M  \begin{array}{c} {\pi} \\ \stackrel{\leftrightharpoons}{\nabla }     \end{array}  N,h   )\label{contr} \end{equation}
of $N$ onto $M$ is a collection of the above data satisfying 
\bea && \pi\nabla=\mbox{id}_M\nn\\
&& D(h)=\mbox{ad}_{d_N}(h)=\nabla\pi-\mbox{id}_N\nn\\
&& \pi h=0,\; h\nabla =0,\; hh=0.\nn\eea
Another way to describe this structure is to say that $M$ is a {\it strong deformation retract (SDR)} of $N$ (also called {\it Eilenberg-Zilber data}). The properties on the last line are referred to as the {\it annihilation properties} or {\it side conditions}. Note that the first line makes $\pi$ surjective (projection) and $\nabla$, injective (inclusion). The map $h$ is also known as the {\it homotopy operator} between $\nabla\pi$ and $\mbox{id}_N$:
\[ \nabla\pi-\mbox{id}_N =D(h)=d_Nh+ hd_N\]
($D=\mbox{ad}_{d_N}$ is the induced differential on $\mbox{Hom}(N,N)$).\end{defi}
Often {\it filtered contractions} are considered.
\begin{remark} Lambe and Stasheff~\cite{LS} noticed that the side conditions on $h$ are not restrictive: if $\pi h=0$ and $h\nabla =0$ are not satisfied, then we can replace $h$ by $h'=D(h)hD(h)$. Now if $h^2=0$ is not satisfied either, we replace $h'$ by $h''=h'd_Nh'$, which finally gives us an operator $h''$ satisfying the side conditions. 
\end{remark}
\begin{lem}\label{A} Given a contraction (\ref{contr}), we have a (not necessarily direct) sum
\[ N=\mbox{Im}(\nabla)+\mbox{Im}(h)+\mbox{Im}(d_N).\]\end{lem}
\noi {\it Proof.} Each $x\in N$ can be written as
\begin{equation} x=\nabla\pi (x)-hd_N(x)-d_Nh(x).\label{split}\end{equation}
$\Box$
\begin{lem}\label{B} \cite{Re} Given a contraction (\ref{contr}), we have 
\[ \mbox{Im}(\nabla)+\mbox{Im}(h)=\mbox{Im}(\nabla)\oplus\mbox{Im}(h)=\mbox{Ker}(h).\]
\end{lem}
{\it Proof.} We have 
\[\mbox{Im}(\nabla)\subset\mbox{Ker}(h)\;\;\;\mbox{and}\;\;\; \mbox{Im}(h) \subset\mbox{Ker}(h) \]
since $ h\nabla =0$ and $h^2=0$. Conversely, by (\ref{split}), each $x\in \mbox{Ker}(h)$ can be written as
\[ x=\nabla\pi (x)-hd_N(x)\in \mbox{Im}(\nabla)+\mbox{Im}(h). \]
That the sum is direct can be seen as follows: let $x\in \mbox{Im}(\nabla)\cap\mbox{Im}(h)$. Then we have $x=\nabla (y)=h(z)$ for some $y\in M$ and $z\in N$. Rewriting the decomposition of $x$ in $\mbox{Ker}(h)$, we obtain
\bea x&=& \nabla\pi (x)-hd_N(x)\nn\\
&=& \nabla\pi h(z)-hd_N\nabla (y)\nn\\
&=& 0-hd_N\nabla (y)\;\;\;(\pi h=0  )\nn\\
&=& h\nabla d_M(y)\;\;\;(\mbox{$\nabla$ chain map})\nn\\
&=& 0 \;\;\;(h\nabla =0).\nn\eea $\Box$
\begin{cor} For any contraction (\ref{contr}), we have
\[  H(N,h)\cong\mbox{Im}(\nabla)\cong M.\]
\end{cor}
\begin{lem}\label{C} Given any contraction (\ref{contr}), we have
\[ \mbox{Im}(h)+\mbox{Im}(d_N)=\mbox{Im}(h)\oplus\mbox{Im}(d_N).\]
Moreover, if $d_M\equiv 0$, then 
\[ \mbox{Im}(h)\oplus\mbox{Im}(d_N)=\mbox{Ker}(\pi).\]\end{lem}
{\it Proof.} Say $x\in  \mbox{Im}(h)\cap\mbox{Im}(d_N)$.
Then $x=h(y)=d_N(z)$ for some $y,z\in N$, 
so that by (\ref{split}) we obtain
\[ x=\nabla\pi h(y)-hd_Nd_N(z)-d_Nhh(y)=0.\]
It is always true that $\mbox{Im}(h)\subset\mbox{Ker}(\pi)$ as $\pi h=0$. If $d_M\equiv 0$, we further have the result
\[ \pi d_N(x)=-d_M\pi x=0,\]
so that altogether
\[ \mbox{Im}(h)\oplus\mbox{Im}(d_N)\subset \mbox{Ker}(\pi).\]
Conversely, for $x\in\mbox{Ker}(\pi)$, we see that (even without the condition $d_M=0$)
\[  \mbox{Ker}(\pi)\subset \mbox{Im}(h)+\mbox{Im}(d_N)\]
since
\[ x= -hd_N(x)-d_Nh(x)\;\;\;\forall x\in \mbox{Ker}(\pi)\]
due to (\ref{split}). $\Box$
\begin{lem}\label{D} For any contraction (\ref{contr}) with $d_M=0$ we have
\[ \mbox{Im}(\nabla)+\mbox{Im}(d_N)=\mbox{Im}(\nabla)\oplus\mbox{Im}(d_N)=\mbox{Ker}(d_N) .\]\end{lem}
\noi {\it Proof.} The given sum is direct: 
let $x\in \mbox{Im}(\nabla)\cap\mbox{Im}(d_N)$.
Then $x=\nabla(y)=d_N(z)$ and by (\ref{split})
\[ x=\nabla\pi d_N(z)-hd_Nd_N(z)-d_Nh\nabla(y)=\nabla\pi d_N(z)=-\nabla d_M\pi (z)=0.\]
Clearly, we have $\mbox{Im}(d_N)\subset \mbox{Ker}(d_N)$. Also $\mbox{Im}(\nabla)\subset \mbox{Ker}(d_N)$, because 
\[ d_N\nabla(x)=-\nabla d_M(x)=0.\]
Conversely, by (\ref{split}),
\[  \mbox{Ker}(d_N) \subset \mbox{Im}(\nabla)+\mbox{Im}(d_N)\]
since we can write
\[ x=\nabla\pi (x)-d_Nh(x)\]
(no condition on $d_M$) for $x\in \mbox{Ker}(d_N)$. $\Box$
\begin{cor} For any contraction (\ref{contr}) with $d_M=0$ we have
\[ H(N,d_N)\cong\mbox{Im}(\nabla)\cong M.\]\end{cor}
\begin{prop} For any contraction (\ref{contr}) with $d_M=0$ we have
\begin{enumerate}
\item 
\[ N=\mbox{Im}(\nabla)\oplus\mbox{Im}(h)\oplus\mbox{Im}(d_N)\]
where
\bea && \mbox{Im}(\nabla)\oplus\mbox{Im}(h)=\mbox{Ker}(h)\nn\\
&& \mbox{Im}(h)\oplus\mbox{Im}(d_N)=\mbox{Ker}(\pi)\nn\\
&& \mbox{Im}(\nabla)\oplus\mbox{Im}(d_N)=\mbox{Ker}(d_N),\nn\eea

\item 
$\mbox{Im}(h)\stackrel{d_N}{\rightarrow} \mbox{Im}(d_N)$ is an isomorphism with inverse $\mbox{Im}(d_N)\stackrel{-h}{\rightarrow} \mbox{Im}(h)$, and

\item
$\mbox{Im}(\nabla)\stackrel{\pi}{\rightarrow}M$ is an isomorphism with inverse $M\stackrel{\nabla}{\rightarrow}\mbox{Im}(\nabla)$; we also have
\[ \mbox{Im}(\nabla)\cong M=\mbox{Im}(\pi) \cong H(N,d_N)\cong H(N,h). \]
\end{enumerate}\end{prop}
\begin{remark} This is a Hodge-type decomposition reminiscent of the case of a compact orientable Riemannian manifold $M$ without boundary. If 
\[ \ast :\Omega^r(M)\rightarrow \Omega^{{\rm dim}(M)-r}\]
is the ``Hodge star operator'' (an isomorphism) and
\[ d:\Omega^{r-1}(M)\rightarrow\Omega^r(M)\]
is the de Rham differential, then we define a ``partial inverse'' $d^{\dagger}$ (the adjoint exterior derivative operator) to $-d$ by $d^{\dagger}=\pm \ast d\ast$. The commutator of $d$ and $d^{\dagger}$ is called the ``Laplace-Beltrami operator'': $\Delta=dd^{\dagger}+d^{\dagger}d$. Then there exists a unique decomposition of the algebra of de Rham forms as follows:
\[ \Omega^r(M)=\mbox{Harm}^r(M)\oplus d(\Omega^{r-1}(M))\oplus d^{\dagger}(\Omega^{r+1}(M)) ,\]
where the ``harmonic forms'' are given by $\mbox{Harm}^r(M)=\mbox{Ker}(\Delta)$. In the case of our general contraction with $d_M=0$, the operators $h$ and $d_N$ replace $d$ and $d^{\dagger}$ respectively. What do we know about $\Delta$ here? We have
\[ \Delta=D(h)=hd_N+d_Nh=(h+d_N)^2=\nabla\pi -\mbox{id}_N.\]
The kernel of this operator is equal to $\nabla(M)$, as we have
\[ (\nabla\pi -\mbox{id}_N)(x)=0\,\Leftrightarrow \, \nabla\pi(x)=x\,\Leftrightarrow \, x\in\nabla(M),\]
or $\mbox{Ker}(\Delta)=\mbox{Im}(\nabla)$. So is there an analog of the Hodge star operator? If we define an isomorphism $ \ast=h+d_N+\mbox{Id}_{\nabla(M)}$ (where the last operator is zero on the remaining direct summands), then we have $\ast^{-1}=-h-d_N+\mbox{Id}_{\nabla(M)}$, and 
\[ \ast d_N\ast^{-1}=(h+d_N+\mbox{Id}_{\nabla(M)})d_N(-h-d_N+\mbox{Id}_{\nabla(M)})=-hd_Nh=-h\mbox{Id}_{{\rm Im}(d_N)}=h.\]
\end{remark}
\begin{remark} The operator $d^{\dagger}$ is more like the BV operator than the (even) Laplacian, which is not square-zero. Another similar case is $Q$ (BRST operator) and $b_0$ (anti-ghost operator), for which we have $Qb_0+b_0Q=L_0$ (the degree operator which is zero on the cohomology).
\end{remark}
\noi {\it Proof.} We only need to prove (2) and part of (3). First, we want to show that $-d_Nh=\mbox{id}_{{\rm Im}(d_N)}$ But then for $x=d_N(y)$, we have
\bea -d_Nh(x)&=&[-d_Nh]d_N(y)\nn\\
&=&[hd_N-\nabla\pi +\mbox{id}_N]d_N(y)\nn\\
&=&-\nabla[\pi d_N](y)+d_N(y)\nn\\
&=&-\nabla[-d_M\pi]d_N(y)+x\nn\\
&=&x.\nn\eea
Similarly, we would like to have $-hd_N=\mbox{id}_{{\rm Im}(h)}$. If $x=h(y)$, then
\bea -hd_N(x)&=& [-hd_N]h(y)\nn\\
&=& [d_Nh-\nabla\pi +\mbox{id}_N]h(y)\nn\\
&=& -\nabla\pi h(y)+h(y)\nn\\
&=& x.\nn\eea
Finally, we have $\pi\nabla =\mbox{id}_M$ and $\nabla\pi\nabla =\nabla\mbox{id}_M=\nabla$, which shows the isomorphism between $\nabla(M)$ and $\pi(N)$. $\Box$ 
\begin{example} Let $(g,d)$ be a chain complex. Assume that the underlying ring is a field. Then there exists a contraction
\begin{equation} (   H(g)  \begin{array}{c} {\pi} \\ \stackrel{\leftrightharpoons}{\nabla } \end{array}g,h) \label{ghomg}\end{equation}
of chain complexes, where the differential on $H(g)$ is zero: we can write $g$ as a linear sum
\[ g=G\oplus \mbox{Ker}(d)=G\oplus \mbox{Im}(d)\oplus H(g)\]
by choosing arbitrary representatives of the homology classes etc.; let us show the decomposition of an element $x$ of $g$ by
\[ x=x_G+x_{{\rm Im}(d)}+x_{H(g)}.\]
Then $\pi$ is the projection of $g$ onto $H(g)$ and $\nabla$ is the inclusion map of $H(g)$ into $g$. Note that as vector spaces $G$ and $\mbox{Im}(d)$ are isomorphic via $d$: let $x,y\in G$. Then
\[ dx=dy\Rightarrow d(x-y)=0\Rightarrow x-y\in \mbox{Ker}(d)\cap G=\{ 0\} \]
and $d:G\rightarrow \mbox{Im}(d)$ is one-to-one as well as onto. We define $h$ to be the inverse of $-d$ on $\mbox{Im}(d)$ and zero on the rest of $g$. The linear map $h$ is square-zero and increases degree by one. Moreover,
\bea &&(dh+hd)(x)\nn\\
&=& (dh+hd)(x_G+x_{{\rm Im}(d)}+x_{H(g)})\nn\\
&=& dh(x_{{\rm Im}(d)})+hd(x_G)\nn\\
&=& -x_{{\rm Im}(d)}-x_G\nn\eea
and
\bea &&(\nabla \pi-\mbox{id}_g)(x_G+x_{{\rm Im}(d)}+x_{H(g)})\nn\\
&=& x_{H(g)}-(x_G+x_{{\rm Im}(d)}+x_{H(g)})\nn\\
&=& -x_{{\rm Im}(d)}-x_G.\nn\eea
In comparison with the last corollary, we have
\bea && (N,d_N)=(g,d)\nn\\
&& (M,d_M)=(H(g,d),0)\nn\\
&& \mbox{Im}(\nabla)=H(g,d)\nn\\
&& \mbox{Im}(h)=G\nn\\
&&\mbox{Im}(d_N)=\mbox{Im}(d).\nn\eea
$\diamond$\end{example}
\begin{example} (The {\it Tensor Trick}) Any contraction~(\ref{contr}) of chain complexes induces a filtered contraction
\[ (   T^c[M]  \begin{array}{c} {T^c\pi} \\ \stackrel{\leftrightharpoons}{T^c\nabla}     \end{array}  T^c[N],T^ch   ) \]
of coaugmented differential graded coalgebras. Here is how: the projection 
\[ P_N:T^c[N]\rightarrow N\]
followed by the surjective chain map $\pi:N\rightarrow M$ gives us a linear map
\bea  &&\pi\circ P_N:T^c[N]\rightarrow M\nn\\
&& (\pi\circ P_N)(x_1\otimes\cdots\otimes x_k)=\left\{ \begin{array}{ll} \pi(x_1) & \mbox{if $k=1$}\\ 0 & \mbox{otherwise.}
\end{array}   \right.\nn\eea
which can then be made into a unique coalgebra map
\[ T^c\pi:T^c[N]\rightarrow T^c[M]\]
with the usual formula
\[ T^c\pi(x_1\otimes\cdots\otimes x_k)=\sum_{i=1}^k x_1\otimes\cdots\otimes \pi(x_i) \otimes\cdots\otimes x_k.\]
Next, the morphisms $T^c\pi$ and $T^c\nabla$ pass to the corresponding morphisms on the coalgebras $S^c[N]$ and $S^c[M]$ respectively, and $S^ch$ is obtained from $T^ch$  by symmetrization, to yield a contraction
\[ (   S^c[M]  \begin{array}{c} {S^c\pi} \\ \stackrel{\leftrightharpoons}{S^c\nabla}     \end{array}  S^c[N],S^ch   ). \]
In particular, the contraction~(\ref{ghomg}) induces 
\begin{equation} (   S^c[sH(g)]  \begin{array}{c} {S^c\pi} \\ \stackrel{\leftrightharpoons}{S^c\nabla } \end{array}S^c[sg],S^ch), \label{sghomg}\end{equation}
which is a filtered contraction of coaugmented DG coalgebras. (Warning: $S^c\pi$ and $S^c\nabla$ are morphisms of coalgebras but $S^ch$ is not a coalgebra morphism, although it is somewhat compatible with the coalgebra structure, being a homotopy of coalgebra maps. One has to be careful when defining a homotopy of cocommutative coalgebras.) $\diamond$
\end{example}

\subsection{The first main theorem.}

Assume that $\partial$ is the codifferential corresponding to an sh-Lie algebra structure on $(g,d)$. Since the corresponding multilinear map on $g$ has other components than the binary bracket, we will denote the symmetric coalgebra on $sg$ with codifferential $\partial$ by $S^c_{\partial}[sg]$ and not by $S^c_{[\, ,\,]}[sg]$. Given two sh-Lie algebras $(g_1,\partial_1)$ and $(g_1,\partial_2)$, an {\it sh-morphism} or {\it sh-Lie map} from $g_1$ to $g_2$ is a morphism $S^c_{\partial_1}[sg_1]\rightarrow S^c_{\partial_2}[sg_2]$ of DG coalgebras.
\vsp

\begin{thm} \label{firstt} Given a DGLA $g$ and a contraction of chain complexes such as (\ref{ghomg}), the data determine

\noi (i) a differential $D$ on $S^c[sH(g)]$ (a coalgebra perturbation of the zero differential) turning the latter into a coaugmented DG coalgebra, hence endowing $H(g)$ with an sh-Lie algebra structure, 

\noi (ii) a Lie algebra twisting cochain 
\[ \tau: S_D^c[sH(g)]\rightarrow g\]
with adjoint $\bar{\tau}$, written
\[ \bar{\tau}=(S^c\nabla)_{\partial}:S^c_D[sH(g)]\rightarrow C[g], \]
that induces an isomorphism on the homology, and

\noi (iii) an extension of $(S^c\nabla)_{\partial}$ to a new contraction
\[ (   S_D^c[sH(g)]  \begin{array}{c} {(S^c\pi)_{\partial}} \\ \stackrel{\leftrightharpoons}{(S^c\nabla)_{\partial}} \end{array}S_{\partial}^c[sg],(S^cd)_{\partial}) \]
of filtered chain complexes (not necessarily of coalgebras). 
\end{thm}

\noi {\bf Notes on Notation.} While the induced bracket on $H(g)$ is a strict graded Lie bracket, the differential $D$ may involve meaningful terms of higher order. Let us introduce a table for the notation used in~\cite{HS} for different types of chain complexes and the corresponding symmetric coalgebras. \vsp

\hspace{-1in}
\[ \begin{array}{|l|l|l|l|l|}
\hline \mbox{Chain complex}&\mbox{Bracket(s)} & \mbox{Sym.~coalgebra} &\mbox{Coderivation} &\mbox{Property}\\ \hline 
{} & {} & {} & {} & {} \\
(g,d)& [\;\; ,\;\; ] & \left( S^c[sg],d    \right)  & \partial & {}\\
\mbox{graded} & \mbox{generic} & \mbox{DG coalgebra; } & \mbox{coderivation} & {}\\
\mbox{chain complex} & \mbox{bilinear bracket}& \mbox{induced diff.~$d$} & \mbox{induced by $[\; ,\; ]$} & {}\\
{} & {} & {} & {} & {} \\ \hline
{} & {} & {} & {} & {} \\
(g,d) & [\;\; ,\;\; ] & \left( S^c_{[\, ,\, ]}[sg]  ,d \right) =C[g] & \partial & \left( d+\partial \right)^2=0\\
\mbox{DGLA} & \mbox{Lie bracket on $g$;} & \mbox{generalized Koszul or} & \mbox{coderivation} & {}\\ 
{} & \mbox{$d$ derivation of it} & \mbox{Chevalley-Eilenberg} & \mbox{induced by $[\; ,\; ]$} & {}\\
{} & {} & \mbox{complex} {} & {} & {} \\
{} & {} & {} & {} & {} \\ \hline
{} & {} & {} & {} & {} \\
(g,0) & [\;\; ,\;\; ] & \left( S^c[sg]  ,0 \right) =\Lambda^{\bullet}(g) & \partial & \partial^2=0\\
\mbox{Lie algebra} & \mbox{Lie bracket} & \mbox{Koszul complex} & \mbox{codifferential} & {}\\
{} & {} & \mbox{for homology} & \mbox{induced by $[\; ,\; ]$} & {}\\
{} & {} & {} & {} & {} \\ \hline
{} & {} & {} & {} & {} \\
(g,d) & \ell_2,\; \ell_3,\dots & \left( S^c_{\partial}[sg]  ,d \right) & \partial & (d+\partial)^2=0\\
\mbox{$L_{\infty}$ algebra;} & \mbox{higher brackets} & {} & \mbox{codiff.~induced by} & {}\\
d=\ell_1 & {} & {} & \ell_2,\; \ell_3,\dots & {}\\
{} & {} & {} & {} & {} \\ \hline
{} & {} & {} & {} & {} \\
(H(g),0) &  [\;\; ,\;\; ] & \left( S^c_{D}[sH(g)],0\right) & D & D^2=0 \\
\mbox{homology of} & \mbox{induced Lie} & {} & \mbox{codiff.~defined} & {}\\
\mbox{DGLA $(g,d)$} & \mbox{bracket on $H(g)$} & {} & \mbox{in the proof} & {}\\
\mbox{with given}& {} & {} & {} & {}\\
\mbox{contraction}{} & {} & {} & {} & {} \\
{} & {} & {} & {} & {} 
\\ \hline \end{array}\]
\pagebreak

\noi {\it Sketch of Proof.} We obtain the differential $D$ and the twisting cochain $\tau$ on $S^c[sH(g)]$ as infinite series by induction: for $b\geq 1$, write $S^c_b$ for the homogeneous degree-$b$ component of $S^c[sH(g)]$. Then $D,\tau$ for $b\geq 2$ are given by
\bea && \tau =\tau^1+\tau^2+\cdots,\;\;\; \tau^1=\nabla{\tau_{H(g)}},\;\;\; \tau^j:S^c_j\rightarrow g,\; j\geq 1,\label{first}\\
&& \tau^b=\frac{1}{2}h([\tau^1,\tau^{b-1}]+\cdots +[\tau^{b-1},\tau^1])\nn\\
&& D=D^1+D^2+\cdots \label{second}\eea
where $D^{b-1}$ is the coderivation of $S^c[sH(g)]$ determined by
\[ \tau_{H(g)}D^{b-1}=\frac{1}{2}\pi ([\tau^1,\tau^{b-1}]+\cdots +[\tau^{b-1},\tau^1]):S^c_b\rightarrow H(g).\]
That is, the coderivation followed by projection onto the degree-one subspace $sH(g)$ of $S^c[sH(g)]$ is given by the above formula. In the notation of Subsection~\ref{twothree}, we have $f_b=\tau_{H(g)}D^{b-1}$ and $D=\hat{f}$.
For example (dropping the symbol $s$ for elements of $sH(g)$), we have $\tau^1(x)=x\in H(g)\subset g$, and $\tau^2(xy)= h[x,y]$ by~(\ref{cupbra}). Let us also compute two terms of $D$:
\[ \tau_{H(g)}D^{1}(xy)=\pi [x,y]\in H(g),\]
and
\[ \tau_{H(g)}D^{2}(xyz)=\frac{1}{2}\pi (\, [x,h[y,z]]+[h[x,y],z]\, ),\]
etc.
We can see why $\tau$ is a LA twisting cochain: since $\tau$ satisfies
\[ \tau =h\left( \frac{1}{2}[\tau,\tau]  \right),\]
we obtain
\[ -d\tau=\frac{1}{2}[\tau,\tau]\]
in case of the particular SDR we constructed, and the last equation is the master equation (the differential on $H(g)$ being zero).
The sums (\ref{first}) and (\ref{second}) are infinite, but when either one is applied to a specific element in some subspace of finite filtration degree, only finitely many terms will be nonzero. (The summand $D^1$ is the ordinary Cartan-Chevalley-Eilenberg differential for the classifying coalgebra of the graded Lie algebra $H(g)$.) The proof that $D$ is indeed a coalgebra differential and $\tau$ is a twisting cochain ``will be given elsewhere''. A ``spectral sequence argument'' shows that $\bar{\tau}$ induces an isomorphism on the homology.$\Box$

\begin{remark}
If $\nabla H(g)$ happens to be a Lie subalgebra of $g$, then $[\tau^1,\tau^1]$ will have values in $H(g)$ and $\tau^2=(1/2)h[\tau^1,\tau^1]$, as well as the remaining $\tau^j$, will be zero. Similarly, we will have $D=D^1$.
\end{remark}

\begin{cor} Under the hypotheses of Theorem~1, 
\begin{equation} \tau: S^c_D[sH(g)]\rightarrow g,\label{mostgen}\end{equation}
viewed as an element of degree $-1$ of the DGLA $\mbox{Hom}(S_D^c[sH(g)],g)$, satisfies the master equation~(\ref{two}).
\end{cor}

The twisting cochain~(\ref{mostgen}) is our most general solution of the master equation. The other solutions of the master equation can be derived from it.

\section{Corollaries and the second main theorem}

\subsection{Other corollaries of Theorem~\ref{firstt}.}

\begin{cor} Under the hypotheses of Theorem~\ref{firstt}, suppose in addition that there is a differential $\tilde{D}$ on $S^c[sH(g)]$ turning the latter into a coaugmented DG coalgebra in such a way that $(S^c\pi){\partial}=\tilde{D}(S^c\pi)$. Then $D=\tilde{D}$ and $(S^c\pi)_{\partial}$ may be taken to be $S^c\pi$. In particular, when $(S^c\pi){\partial}$ is zero, then the differential $D$ on $S^c[sH(g)]$ is necessarily zero, that is, the new contraction in Theorem~1 has the form 
\[ (   S^c[sH(g)]  \begin{array}{c} {S^c\pi} \\ \stackrel{\leftrightharpoons}{(S^c\nabla)_{\partial}} \end{array}S_{\partial}^c[sg],(S^ch)_{\partial}) .\]
For example, this is the case when the composite 
\[ g\otimes g \stackrel{[\, ,\,]}{\rightarrow}g \stackrel{\pi}{\rightarrow} H(g) \]
is zero.
\end{cor}

\begin{cor} Under the hypotheses of Theorem~\ref{firstt}, suppose in addition that there is a differential $\tilde{D}$ on $S^c[sH(g)]$ turning the latter into a coaugmented DG coalgebra in such a way that $\partial(S^c\nabla)=(S^c\nabla)\tilde{D}$. Then $D=\tilde{D}$ and $(S^c\nabla)_{\partial}=S^\nabla$. In particular, when $\partial(S^c\nabla)$ is zero, then the differential $D$ on $S^c[sH(g)]$ is necessarily zero, that is, the new contraction in Theorem~1 has the form 
\[ (   S^c[sH(g)]  \begin{array}{c} ({S^c\pi})_{\partial} \\ \stackrel{\leftrightharpoons}{S^c\nabla} \end{array}S_{\partial}^c[sg],(S^ch)_{\partial}) .\]
For example, this is the case when the composite 
\[ H(g)\otimes H(g)\rightarrow \stackrel {\nabla\otimes\nabla}{\rightarrow} g\otimes g\stackrel{[\, ,\,]}{\rightarrow}g\]
is zero.
\end{cor}

\subsection{The second main theorem}

\begin{thm}\label{secondd} Given a DGLA $g$, a DGL subalgebra $m$ of $g$, and a contraction
\[ (   H(g)  \begin{array}{c} {\pi} \\ \stackrel{\leftrightharpoons}{\nabla } \end{array}m,h) \]
of chain complexes so that the composite
\[ m\otimes m \stackrel{[\, ,\,]}{\rightarrow}m \stackrel{\pi}{\rightarrow} H(g) \]
is zero, then the induced bracket on $H(g)$ is zero, that is, $H(g)$ is abelian as a graded Lie algebra, and the data determine a solution $\tau\in\mbox{Hom}(S^c[sH(g)],g)$ of the master equation~(\ref{two}) in such a way that the following hold:

\noi (i) The composite $\pi\tau$ coincides with the universal twisting cochain $S^c[sH(g)]\rightarrow H(g)$ for the abelian Lie algebra $H(g)$, and

\noi (ii) the values of $\tau$ lie in $m$.
\end{thm}

%\section{Ideas}
%
%Use formulas involving $\Phi$, $\Delta$, and $Q$ from old BV paper to prove our result about sh-Lie on the abelian cohomology in a different way. How can we split $V$ so that $H$ is directly comm assoc? We are already given bilinear and trilinear maps that make $m$ homotopy comm and assoc. Im($Q$) is fixed, but can we finagle the rest? In fact, we wanna show either the resulting sh-Lie is the same as the $\Phi$'s or ignore $\Phi$'s altogether. Or: start from the opposite direction. Look at the coderivation corresponding to $\Phi^2+\cdots$.

\section{Differential Gerstenhaber and BV algebras}

\subsection{Differential Gerstenhaber algebras}

\begin{defi} A {\it Gerstenhaber (or G-) algebra} consists of \begin{itemize}
\item A graded commutative and associative algebra $(A,\mu)$ ($\mu$ suppressed), and
\item A graded Lie bracket (the {\it Gerstenhaber} or $G-$ {\it bracket}) $[\, ,\, ]:A\otimes A\rightarrow A$ of degree $-1$, such that
\item For each homogeneous element $a\in A$, bracketing with $a$ is a derivation of the Lie bracket of degree $|a|-1$.
\end{itemize}

That is, we want $\mbox{ad}(a)$ to commute with the bracket for all $a\in A$.\end{defi}

\begin{defi} A {\it differential G-algebra} is a Gerstenhaber algebra $(A,[\, ,\, ])$ with a differential $d$ of degree $+1$ on $A$ which is a derivation of the multiplication on $A$. 
\end{defi}

We want $[d,\mu]=0$ and $[d,d]=0$ in Gerstenhaber's composition bracket notation.

\vsp

\begin{defi} A differential G-algebra is called {\it strict} if the differential $d$ is a derivation of the $G$-bracket as well.\end{defi}

We want $d$ to commute with the bracket.\vsp

Let $(A,[\, ,\, ],d)$ be a strict differential G-algebra. We will for the moment ignore all the extra structure on $A$ except for the G-bracket and the differential. As such, $A$ is a DGLA, and we will change the notation to $g$ to emphasize that. We will use the grading
\[ g_1=A^0,\; g_{0}=A^1,\; g_{-1}=A^2, \dots, g_{-n}=A^{n+1},\dots \]
so that the graded bracket and the differential on $g$ are now ``ordinary'': namely,
\[ [\, ,\, ]:g_{j}\otimes g_k\rightarrow g_{j+k},\;\;\; d:g_j\rightarrow g_{j-1}.\]
Consider a contraction of $g$ onto $H(g)$ as in~(\ref{ghomg}). Let $\partial$ denote the operator (``perturbation of $d$'') on $S^c_{\partial}[sg]$ corresponding to the Lie bracket on $g$. By the Main Theorem (Theorem~1), we can transfer it to the symmetric coalgebra of the homology: there exists a new contraction
\[ (   S_D^c[sH(g)]  \begin{array}{c} {(S^c\pi)_{\partial}} \\ \stackrel{\leftrightharpoons}{(S^c\nabla)_{\partial}} \end{array}S_{\partial}^c[sg],(S^cd)_{\partial}) \]
of not only filtered chain complexes but of {\it filtered differential graded coalgebras}. The twisting cochain $\tau$ of~(\ref{mostgen}) of Corollary~1 is now an element of
\[ \mbox{Hom}(S^c_D[sH(g)],s^{-1}g)\]
of degree $-2$, satisfying the master equation
\[ D\tau =\frac{1}{2}[\tau,\tau],\] where $D$ is the Hom-differential and the graded cup bracket on the right-hand side refers to the one induced by the graded coalgebra structure on $S^c[sH(g)]$ and the graded Lie algebra structure on $g$. %(It is also possible to rewrite $\tau$ in terms of cohomology degrees, but we will skip that.)

\subsection{Differential BV algebras}

\begin{defi} Let $(A,[\, ,\, ])$ be a G-algebra with an additional operator $\Delta$ on $A$ of degree $-1$. If $\Delta$ satisfies the condition
\[ [a,b]=(-1)^{|a|}\left(  \Delta(ab)-(\Delta a)b-(-1)^{|a|}a(\Delta b)\right) ,\]
then it is said to be a {\it generator} of the G-algebra. In this case, $(A,\Delta)$ is called a {\it weak Batalin-Vilkovisky (BV-) algebra}. If, moreover, $\Delta$ is exact (i.e.~$\Delta^2=0$), then $(A,\Delta)$ is simply called a {\it Batalin-Vilkovisky (BV-) algebra}. Koszul has shown that $\Delta$ behaves as a derivation for the G-bracket:
\[ \Delta [x,y]=[\Delta x,y]-(-1)^{|x|}[x,\Delta y]\;\;\; \forall\, x,y\in A.\]
With respect to the original graded commutative and associative product on $A$, we can only say that $\Delta$ is a second order differential operator, or $\Phi^{3}_{\Delta}(a,b)=0$, where $\Phi^r_{\Delta}$ are $r$-linear operators used to define higher order differential operators.
\end{defi} 

Now let us denote the bracket on a (weak) BV-algebra $(A,\Delta)$ by $[\, ,\, ]$. 

\begin{defi} If $d$ is a differential of degree $+1$ that endows $(A,[\, ,\, ])$ with a differential G-algebra structure, such that $\Delta d+d\Delta =0$, then the triple $(A,\Delta,d)$ is called a {\it (weak) differential BV-algebra}.
\end{defi}

\begin{prop} For any weak differential BV-algebra $(A,\Delta,d)$, the differential $d$ behaves as a derivation of the G-bracket $[\, ,\, ]$: 
\[ d[x,y]=[dx,y]-(-1)^{|x|}[x,dy].\]
That is, $(A,[\, ,\, ],d)$ is a differential G-algebra.
\end{prop}

Note: $(A,[\, ,\, ],\Delta)$ is not a differential G-algebra unless $\Delta$ is exact.
\vsp

\begin{example} \cite{Fi} Let $V$ be a $\mbox{\bf Z}_2$-graded finite dimensional vector space and $sV^{\ast}$ be its suspended-graded dual. If $\{ x_i\}$ is a basis for $V$ consisting of homogeneous elements, then the dual basis $\{ x_i^{\ast}\}$ has the property that $x_i$ and $x_i^{\ast}$ always have opposite parities. Then the algebra $\mbox{C}[[x_1,\dots,x_n,x_1^{\ast},\dots,x_n^{\ast}]]$ of formal power series has the following BV operator (Laplacian):
\[ \Delta =\sum_{i=1}^{n}\frac{\partial}{\partial x_i^{\ast}}\,\frac{\partial}{\partial x_i}.\]
Since the underlying algebra is graded commutative, the composition of two derivations is a second order differential operator by any definition. Moreover $\Delta^2=0$, which makes a BV-algebra out of this data. $\Diamond$ \end{example}

%\noi {\bf Example.} Let $(A,m)$ be a finite dimensional associative algebra, and $\{ a_i\}$ be a basis of $A$ with dual basis $\{ a_i^{\ast}\}$ of $A^{\ast}$. Then the multiplication $m$ can be written as
%\[ m=\sum_{i,j,k}\,\lambda_{ijk}\, a_i\,\otimes a_j^{\ast}\otimes a_k^{\ast},\]
%and the (derived) Lie bracket 
%\[ [a,b]_m=[[m,a],b]=m(a,b)-m(b,a)\]
%is the BV-bracket in the Hochschild complex (where elements of $A$ are also considered multilinear operators) obtained by bracketing $m$ with first $a$ then $b$. %Clearly, $m$ has the form
%\[ m=\sum_{i,j,k}\lambda_{ijk}\, a_i\,\frac{\partial}{\partial a_j}\,\frac{\partial}{\partial a_k}.\]
%\vsp

%\begin{example} Khudaverdian and Nersessian say that the BV operator on a super K\"{a}hlerian manifold corresponds to the {\it covariant divergence} operator $\ast d\ast$.
%\end{example}

\subsection{Formality}

%\subsubsection{Homotopy}
%
%\noi {\bf Definition.} Let $C,D$ be chain complexes and let $F,G:C\rightarrow D$ be chain maps (i.e. we have $F_n:C_n\rightarrow D_n$, same with $G_n$, such that $Fd=dF$ and $Gd=dG$). A {\it chain homotopy} $S$  from $F$ to $G$ consists of
%a homomorphism
%\[ S_n:C_n\rightarrow D_{n+1}\]
%for each $n$ such that
%\[ dS+Sd=F-G .\]

\subsubsection{Formality of differential graded $P$-algebras}

Recall that our ground ring is a field of characteristic zero. Let $P$ be a differential graded operad and $(A,d)$ be an algebra over $P$. We often want to know to which extent the cohomology of a space reflects the underlying topological or geometrical properties of that space.\vsp

\begin{defi} The $P$-algebra $A$ is called {\it formal} if there exists a strongly homotopy $P$-algebra map $(H(A),0)\stackrel{F}{\rightarrow} (A,d)$ which induces an isomorphism in homology. 
\end{defi}

%I asked Stasheff why this property is called ``formality''. He said ``not sure it was that thoughtful... something like
%{\it the homotopy type is a formal consequence of the cohomology}''.
%\vsp

%Following Tamarkin and Tsygan: (1)~An object from classical calculus should have a noncommutative analogue (for example, we have $A=C^{\infty}(M)$ with noncommutative analogue $\mbox{Der}(A)=\mbox{Vect}(M)$). (2)~Any algebra structure on the classical object should be present in the noncommutative version up to strong homotopy. [This should be continued...]

\subsubsection{Examples}

\begin{example} {\bf (Commutative DG  associative algebras.)} A smooth manifold $M$ is called {\it formal} if the commutative associative DG algebra of de Rham forms on $M$ is formal in the sense of the above Definition. Examples are compact K\"{a}hler manifolds, Lie groups, and complete intersections. Poisson manifolds (proof by Sharygin and Talalaev). $\diamond$
\end{example}

\begin{example} {\bf (DGLA's.)} The Hochschild complex for the algebra $A=C^{\infty}(M)$ of smooth functions on a Poisson manifold $M$ (Kontsevich). $\diamond$
\end{example}

\subsection{Differential BV algebras and formality}

\begin{defi} We will say a differential BV-algebra $(A,\Delta,d)$ ($A\! =\! g$ as a Lie algebra) {\it  satisfies the statement of the K\"{a}hlerian Formality Lemma} (or the $\partial\bar{\partial}$ Lemma) if the maps
\[ \left( \mbox{Ker}(\Delta),d\, |_{\mbox{\small Ker}(\Delta)}  \right)\hookrightarrow(g,d),\;\;\; \left( \mbox{Ker}(\Delta),d\, |_{\mbox{\small Ker}(\Delta)}    \right)\stackrel{\mbox{\small proj}}{\longrightarrow} H(g,\Delta)\]
are isomorphisms on the homology, where $H(g,\Delta)$ is endowed with the zero differential. \end{defi}

\begin{remark} If the statement of the K.F.L.~is satisfied, then proj can be extended to a contraction
\[  \left(  ( H(g,d),0)  \begin{array}{c} {\pi} \\ \stackrel{\leftrightharpoons}{\nabla } \end{array}\left( \mbox{Ker}(\Delta),d\, |_{\mbox{\small Ker}(\Delta)}\right),h\right) .  \]
Since we have
\[ H\left( \mbox{Ker}(\Delta),d\right)\cong H(g,d)\]
and
\[ H\left( \mbox{Ker}(\Delta),d\right)\cong H(g,\Delta),\]
now we can have a contraction of $\mbox{Ker}(\Delta)$ onto $H(\mbox{Ker}(\Delta))=H(g)$ as we did with $g$ and $H(g)$. 
\end{remark}

\begin{thm} \label{thirdd} Let  $(A,\Delta,d)$ be a differential BV-algebra satisfying the statement of the K\"{a}hlerian Formality Lemma and extend the projection proj to a contraction
\[  \left(   (H(g,d),0)  \begin{array}{c} {\pi} \\ \stackrel{\leftrightharpoons}{\nabla } \end{array}m,h\right)   \]
where
\[ m=\left( \mbox{Ker}(\Delta),d\, |_{\mbox{\small Ker}(\Delta)}\right) \]
and $\pi =\,${\em proj}. Then $H(g)$ is \underline{abelian} as a graded Lie algebra and the data determine a solution 
\[ \tau\in\mbox{Hom}(S^c[sH(g)],g) \]
of the master equation $\displaystyle{d\tau =\frac{1}{2}[\tau,\tau]}$ in such a way that the following hold:
\begin{itemize}
\item The values of $\tau$ lie in $m$, that is, the composite
\[ \Delta\circ\tau :S^c[sH(g)]\rightarrow g\]
is zero;
\item The composite $\pi\tau$ coincides with the universal twisting cochain for the abelian graded Lie algebra $H(g)$; so that
\item For $k\geq 2$, the values of the component $\tau_k$ of $\tau$ on $S^c_k[sH(g)]$ lie in Im$(\Delta)$.
\end{itemize}\end{thm}

\noi {\it Proof.} Follows from Theorem~\ref{secondd}. $\Box$
\vsp

Let us add the following condition to the ones in Theorem~\ref{thirdd}: suppose that $A$ consists of a single copy of the ground ring $F$ (necessarily generated by the unit 1 of $A$) and that $\Delta(1)=0$. Then 1 generates a central copy of the ground ring in $g$ (the ground ring commutes with all elements of $A$), and we may write
\[ g=F\oplus \tilde{g}\]
as a direct sum of differential graded Lie algebras. Here $\tilde{g}$ is the uniquely determined complement of $F$. Why unique? We have $g_{1}=A_0=F$ and we may take 
\[ \tilde{g}=\bigoplus_{n\geq 0}g_{-n}=\bigoplus_{n\geq 0}A^{n+1},\] 
where $\tilde{g}$ will be closed under the (degree-zero) Lie bracket:
\[ [\, ,\, ]: g_j \otimes g_k \rightarrow g_{j+k} ,\]
with $j+k\leq 0$ if $j,k\leq 0$.\vsp

\begin{cor} Assume that the hypotheses of Theorem~\ref{thirdd}, 
the abovementioned conditions ($A_0=F$, $\Delta(1)=0$), and the condition $H_1(g)\neq 0$ hold. Then the contraction of the Theorem can be chosen in such a way that 
\[ \mbox{Im}(\tau_k)\subset \tilde{g}\;\;\;\mbox{for}\;\;\; k\geq 2.\]
\end{cor}

%\noi {\it Proof.} (Omitted) $\Box$

The statements of the main theorems have an interpretation in the
context of deformation theory, as explained below.

% **************************************************************
\section{Deformation theory}
% Solutions provided by HS Theorems
Given a DGLA $g$, the construction of the universal solution $\tau$ from
Theorems \ref{secondd} and \ref{thirdd} relies on a chosen contraction.
This provides a formal solution of the master equation (MCE),
with a perturbed differential $\C{D}$ on $S^c[sHg]$ in the direction of 
(starting with) the Lie bracket induced on homology,
endowing the former with a dg-coalgebra structure and a twisting cochain:
$$\tau:S^c_\C{D}[s Hg]\to g.$$
% Deformation theory interpretation
The moduli space interpretation of the set of solutions is along the lines of Schlesinger-Stasheff \cite{SS}. Since our focus is on the construction of solutions of the MCE, the reader is referred to the original text \cite{HS}. Additional details in terms of deformation functors, tangent cohomology, and the Kuranishi functor can be found in \cite{Manetti}. The relation between the latter functor and the construction of a twisting cochain corresponding to a contraction will be investigated elsewhere.

%\pagebreak

% **********************************************************************
%		bibliography
% **********************************************************************

\end{document}